\documentclass[12pt]{article}
\usepackage{amsmath,amsthm,amscd,amssymb}
\usepackage{latexsym,epsf,epsfig}
\usepackage[all]{xy}
\usepackage{graphicx}

\newtheorem{proposition}{Proposition}[section] 
\newtheorem{theorem}[proposition]{Theorem} 
 
\newtheorem{corollary}[proposition]{Corollary}

\title{On the Area of Pedal and Antipedal Triangles}
\author{Adrian Mitrea}        

\date{}    

\begin{document}

\maketitle

\abstract
{We give a new proof of the formula expressing the area of the triangle 
whose vertices are the projections of an arbitrary point in the plane onto 
the sides of a given triangle, in terms of the geometry of  
the given triangle and the location of the projection point. 
Other related geometrical constructions and formulas are also presented.}

\section{Introduction}

The setting in which all results of this paper are stated is that
of a two-dimensional Euclidean plane. Given a point 
$P$ in the plane of a given triangle
$ABC$, denoted by $\Delta ABC$, a triangle $A'B'C'$ is called the {\it pedal triangle
of $P$ with respect to $\Delta ABC$} if $A'$, $B'$ and $C'$
are the projections of $P$ onto the lines $BC$, $AC$, and $AB$,
respectively. The point $P$ will be called the {\it pedal point} of 
$\Delta A'B'C'$ with respect to $\Delta ABC$ (see e.g., \cite{Co}, 
\cite{Jo}). Figures~1 and~2 depict two examples of pedal 
triangles, corresponding
to the point $P$ being inside and outside of $\Delta ABC$, respectively.

\begin{center}
\includegraphics[scale=0.8]{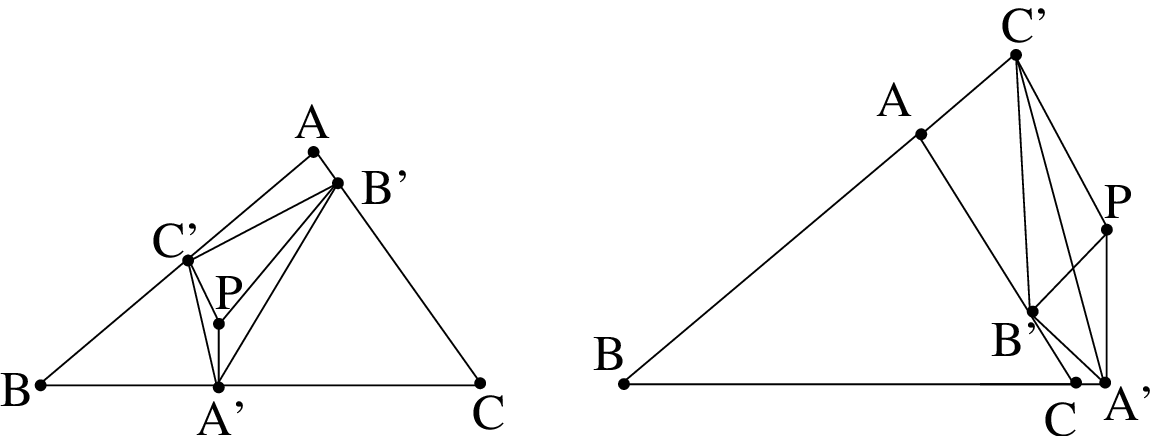}
\end{center}
\centerline{{\bf Figure~1}\hskip 1.8in {\bf Figure~2}}
\vskip 0.15in

\noindent A familiar instance
of a pedal triangle is the {\it orthic triangle} of a given triangle,
the pedal triangle of the {\it orthocenter} of the given triangle, the 
meeting point of its three altitudes.

The main goal of this note is to give a new proof to the formula for the 
area of a pedal triangle of a point, relative to a fixed triangle. 
This formula takes into consideration, besides the geometrical 
characteristics of this fixed triangle, only the location of the pedal point. 
With the convention that $|\Delta XYZ|$ denotes the area
of the triangle $XYZ$, the following holds:

\begin{theorem}\label{luigi}
Let $\Delta ABC$ be a given triangle in the plane, and denote by $O$ and $R$ the center 
and the radius of the circumcircle, respectively. Let $P$ be an arbitrary 
point in the plane of $\Delta ABC$, and let $\Delta A'B'C'$ be the pedal 
triangle of $P$ with respect to $\Delta ABC$. Then

\begin{eqnarray}\label{supermario}
\frac{|\Delta A'B'C'|}{|\Delta ABC|}=\frac{|R^2-OP^2|}{4R^2}.
\end{eqnarray}
\end{theorem}  

This is a classical result.
However, the proofs existing in the literature (we are aware of
\cite{Ga} and \cite{Op}) are rather complex and involved. Here we present an 
approach of algebraic nature, which is considerably more economical 
and direct. In addition, as consequences of Theorem~\ref{luigi} 
we note a couple of results of independent interest. 

\begin{corollary}\label{tVC-1}
Given a fixed triangle $ABC$, denote by $\Delta A'B'C'$ the pedal triangle of 
an arbitrary point $P$ with respect to $\Delta ABC$.  Then the locus of all points 
$P$ such that the ratio of the area of $\Delta A'B'C'$ to that of $\Delta ABC$ is 
a fixed, given constant, is a circle concentric with the circumcircle of $\Delta ABC$.
\end{corollary}

\begin{corollary}\label{tVC-2}
The locus of all points with the property that their projections 
onto the sides of a given triangle $ABC$ are three collinear points 
is the circumcircle of $\Delta ABC$. 
\end{corollary}

\noindent These are both obvious from (\ref{supermario}). Corollary~\ref{tVC-2}
is usually attributed to Simson, and our contribution is 
to provide a conceptually new proof of this well-known fact.

$\Delta TUV$ is called the {\it antipedal triangle} of a point $K$ with 
respect to $\Delta\,ABC$ if the lines $KA, KB, KC$ are perpendicular to 
$VU$, $TV$ and $TU$, at points $A$, $B$, and $C$, respectively. 
Consequently, it follows that $\Delta TUV$ is the antipedal
triangle of a point $K$ with respect to $\Delta\,ABC$ if and only if 
$\Delta\,ABC$ is the pedal triangle of the point $K$ with respect to 
$\Delta TUV$. In the examples from Figure~1 and Figure~2, $\Delta ABC$ is the antipedal
triangle of point $P$ with respect to $\Delta A'B'C'$.

For a given triangle $ABC$, line $AD'$ is {\it isogonal} to line $AD$ in
${\angle\!\!\!)}\,BAC$ if the bisector of ${\angle\!\!\!)}\,BAC$ is also the
bisector of ${\angle\!\!\!)}\,DAD'$; that is, if $AD'$ is the reflection of $AD$ 
in the bisector of ${\angle\!\!\!)}\,BAC$.  If lines $AD$, $BE$, and $CF$ meet at
a point $K$ then it is known that their isogonals $AD'$, $BE'$, and $CF'$ in their 
respective angles $BAC$, $BCA$, and $BCA$ concur in a point $K'$.  The point $K'$ is
called the {\it isogonal} of the point $K$ for the triangle $ABC$. As a natural 
counterpart to Theorem~\ref{luigi} we also derive a similar formula for the area of 
an antipedal triangle.

\begin{theorem}\label{tata5}
If $\Delta\,TUV$ is the antipedal triangle of the point $K$ with 
respect to $\Delta\,ABC$ then
\begin{eqnarray}\label{boise3}
\frac{|\Delta TUV|}{|\Delta ABC|}=\frac{4R^2}{|R^2-{OK_1}^2|}
\end{eqnarray}
\noindent with $O$ being the center and
$R$ the radius of the circumcircle of $\Delta ABC$, 
and $K_1$ being the isogonal of $K$ (see Figure~6).
\end{theorem} 

\section{The Proof of Theorem~\ref{luigi}}

The idea of this proof is to use analytic geometry in order to recast 
(\ref{supermario}) as a quadratic equation in the coordinates
of the variable point $P$.  Note that, in this scheme, it is not 
necessary to carefully keep track of the specific values
of the coefficients of the quadratic equation; indeed, it is only 
the very nature of the algebraic equation which plays a role.  This is
a general principle which could be useful for other types of problems
as well.  

Turning to specifics, consider the lines $AB$, $BC$, $AC$,
given by the equations $\alpha_C x+\beta_C y +\gamma_C=0$,
$\alpha_A x+\beta_A y +\gamma_A=0$, $\alpha_B x+\beta_B y +\gamma_B=0$,
respectively, where the signs of the corresponding coefficients for each 
line are selected such that if a point $P(x,y)$ is inside $\Delta ABC$, 
then $\alpha_C x+\beta_C y +\gamma_C>0$,
$\alpha_A x+\beta_A y +\gamma_A>0$, $\alpha_B x+\beta_B y +\gamma_B>0$.  
Also, for a point $P(x_1,y_1)$, we denote by $d_C$, $d_A$, and $d_B$
the {\it directed} distance from $P$ to $AB$, $BC$, and $AC$, respectively.
As a result, we have explicit formulas for $d_C$, $d_A$, and $d_B$.  For
example, 
\begin{eqnarray}\label{mario222}
d_C=\frac{\alpha_C x_1+\beta_C y_1 +\gamma_C}
{\sqrt{\alpha_C^2+\beta^2_C}}, 
\end{eqnarray}
and similar expressions hold for $d_A$ and $d_B$.
We note that the area of $\Delta A'B'C'$ can be written as a
linear combination of the areas of $\Delta A'PB'$, $\Delta A'PC'$, 
and $\Delta C'PB'$.  More specifically, one has 
\begin{eqnarray}\label{FFF9}
|\Delta A'B'C'|=\pm|\Delta A'PB'|\pm|\Delta A'PC'|\pm|\Delta C'PB'|
\end{eqnarray}
\noindent where the selection of $+$ or $-$ is dictated by the location of
the point $P$. For example,
if $P$ is inside $\Delta ABC$ (as it is the case in Figure~1), 
then the area of $\Delta A'B'C'$ is equal to the sum of the areas 
of $\Delta A'PB'$, $\Delta A'PC'$, and $\Delta C'PB'$, which 
further implies that the signs of the each of the terms in
(\ref{FFF9}) should be $+$. Another example, corresponding to the
point $P$ as in Figure~2, leads to the formula
$|\Delta A'B'C'|=|\Delta A'PB'|-|\Delta A'PC'|+|\Delta C'PB'|$,
so the corresponding signs are $(+,-,+)$.
A specific choice of the three signs in the right-hand
side of (\ref{FFF9}) turns out to depend on the location
of the point $P$ relative to the three lines $AB$, $BC$, and $AC$,
as well as the circumcircle of $\Delta ABC$.  
%
%\begin{center}
%\includegraphics[scale=1]{concifigure.eps}
%\end{center}
%\centerline{{\bf Figure~3}}
%\vskip 0.15in
%
Figure~3 shows all
possible combinations of signs associated with the various sub-regions
in which the plane is partitioned by the aforementioned lines and
circle.
\begin{center}
\includegraphics[scale=0.7]{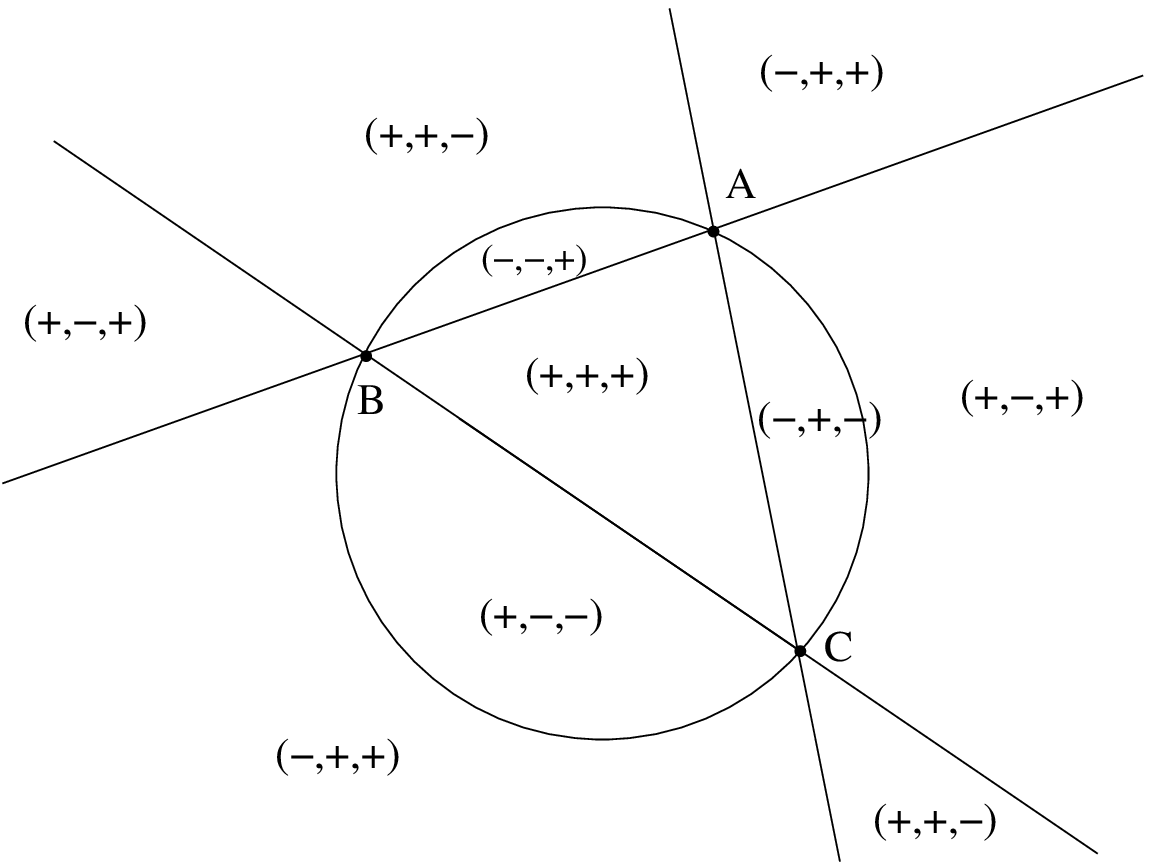}
\end{center}
\centerline{{\bf Figure~3}}
\vskip 0.15in
The directed distances introduced can be used to express the areas of each
of the triangles in the right hand side of (\ref{FFF9}) as follows:
\begin{eqnarray}\label{tf}
&&|\Delta PB'C'|=\frac{|d_B\,d_C\sin({\angle\!\!\!)}\,A)|}{2},\quad
|\Delta A'PB'|=\frac{|d_A\,d_B\sin({\angle\!\!\!)}\,C)|}{2},
\\[4pt]
&&|\Delta A'PC'|=\frac{|d_A\,d_C\sin({\angle\!\!\!)}\,B)|}{2}.
\label{tf2}
\end{eqnarray}
The absolute values from the numerators in these formulas can be dropped
by keeping careful track of the signs of the directed distances based on the
position of $P$. Figure~4 shows the signs of each component in the triplet
$\bigl(d_Bd_C,d_Ad_C,d_Ad_B\bigr)$ associated with the various sub-regions
in which the plane is partitioned by the lines $AB$, $AC$, and $BC$.

\begin{center}
\includegraphics[scale=0.6]{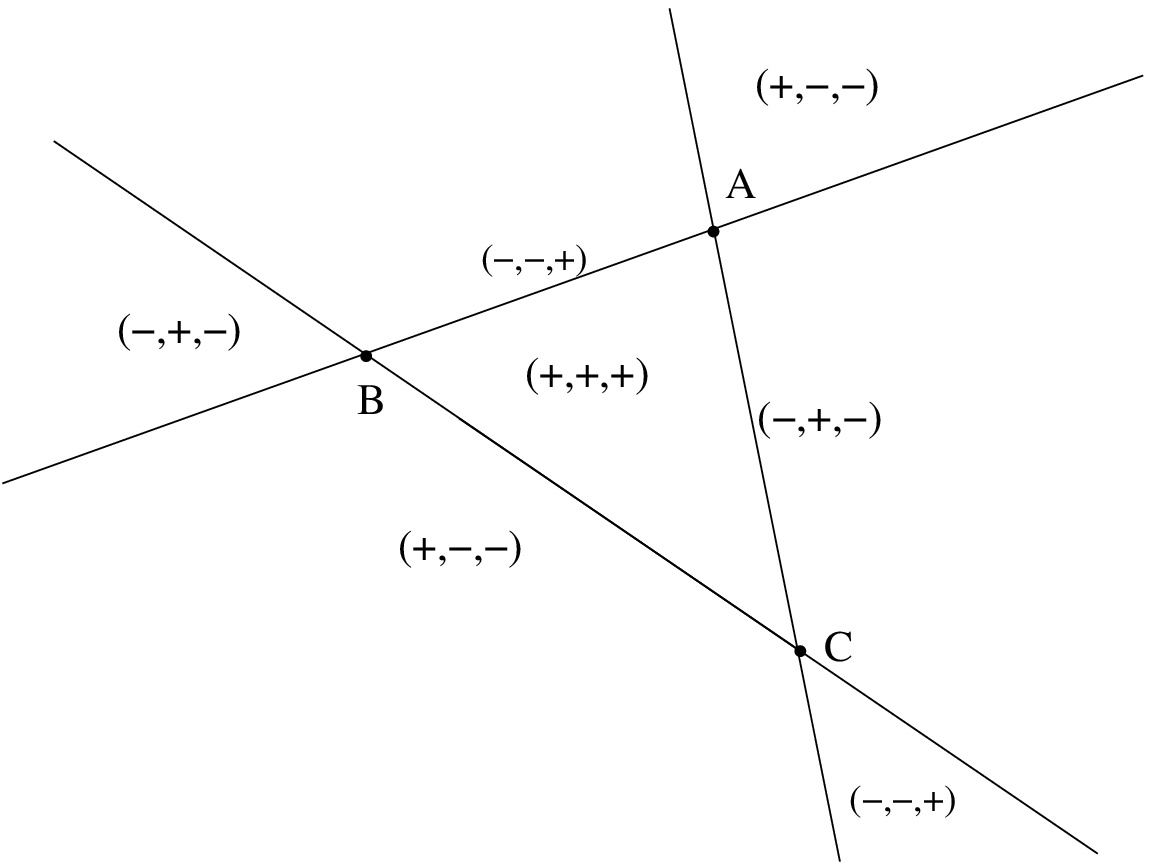}
\end{center}
\centerline{{\bf Figure~4}}
\vskip 0.15in

What is remarkable is the fact that after combining the signs of the 
fractions in (\ref{FFF9}), as depicted by Figure~3,
with the signs of the products of the directed distance, 
as depicted by Figure~4, formula (\ref{FFF9}), in junction with
(\ref{tf})-(\ref{tf2}), becomes
\begin{eqnarray}\label{FFF10}
\pm|\Delta A'B'C'|=\frac{d_B\,d_C\sin({\angle\!\!\!)}\,A)}{2}
+\frac{d_A\,d_B\sin({\angle\!\!\!)}\,C)}{2}
+\frac{d_A\,d_C\sin({\angle\!\!\!)}\,B)}{2},
\end{eqnarray}
\noindent where $+$ corresponds to the case when $P$ is contained in
the circumcircle of $\Delta ABC$, with center $O$ (denoted by $(O)$), and $-$ corresponds to the 
case when $P$ is outside $(O)$.  Making now use of (\ref{mario222}) and
the corresponding formulas for $d_A$, $d_B$, we can
re-write (\ref{FFF10}) as 
\begin{eqnarray}\label{mario1111}
\pm|\Delta A'B'C'| &=& \frac{\alpha_B x_1+\beta_B y_1 
+\gamma_B}{\sqrt{\alpha_B^2+\beta^2_B}}\cdot\frac{\alpha_C x_1+\beta_C y_1 
+\gamma_C}{\sqrt{\alpha_C^2+\beta^2_C}}\cdot\frac{\sin({\angle\!\!\!)}\,A)}{2}
\nonumber\\[4pt]
&& +\frac{\alpha_A x_1+\beta_A y_1+\gamma_A}{\sqrt{\alpha_A^2+\beta^2_A}}
\cdot\frac{\alpha_C x_1+\beta_C y_1+\gamma_C}{\sqrt{\alpha_C^2+\beta^2_C}}
\cdot\frac{\sin({\angle\!\!\!)}\,B)}{2}
\nonumber\\[4pt]
&& +\frac{\alpha_A x_1+\beta_A y_1+\gamma_A}{\sqrt{\alpha_A^2+\beta^2_A}}
\cdot\frac{\alpha_B x_1+\beta_B y_1+\gamma_B}
{\sqrt{\alpha_B^2+\beta^2_B}}\cdot\frac{\sin({\angle\!\!\!)}\,C)}{2}.
\end{eqnarray}
\noindent In addition, using the fact that $|\Delta ABC|$ 
is a real constant that depends only on $A$, $B$ and $C$, 
(\ref{mario1111}) yields
\begin{eqnarray}\label{mario3}
\pm\frac{|\Delta A'B'C'|}{|\Delta ABC|}&=&(ax_1+by_1+c)(dx_1+ey_1+f)
\nonumber
\\
&&
+(gx_1+hy_1+i)(jx_1+ky_1+l)
\nonumber
\\[4pt]
&&
+(mx_1+ny_1+o)(px_1+qy_1+r),
\end{eqnarray}
\noindent where $a,b,c,d,e,f,g,h,i,j,k,l,m,n,o,p,q$ and $r$ are real 
constants that depend only on $A$, $B$ and $C$. 

At this point, we observe that (\ref{supermario}) becomes
\begin{eqnarray}\label{supermariob}
\pm\frac{|\Delta A'B'C'|}{|\Delta ABC|}=\frac{R^2-OP^2}{4R^2},
\end{eqnarray}
\noindent provided we select $+$ when $P$ is in $(O)$ and $-$
when $P$ is outside $(O)$. Clearly, the right hand side of
(\ref{supermariob}) is a quadratic expression in $x_1$ and $y_1$.
Hence, if we now take into account 
(\ref{supermariob}) and (\ref{mario3}), we obtain that 
(\ref{supermario}) is equivalent with 
\begin{eqnarray}\label{mariokart}
\lambda_1x^2_1+\lambda_2y^2_1+\lambda_3x_1y_1+\lambda_4x_1+\lambda_5y_1
+\lambda_6=0,
\end{eqnarray}
\noindent where $\lambda_1$, $\lambda_2$, $\lambda_3$, $\lambda_4$, 
$\lambda_5$, and $\lambda_6$ are real constants that depend only on 
$A$, $B$, and $C$. Any points that satisfy (\ref{supermario}) satisfy 
(\ref{mariokart}), and vice versa. Furthermore, it is fairly easy to 
see that the point $O$ satisfies (\ref{supermario}), as 
both sides of (\ref{supermario}) will be $\frac{1}{4}$.  In addition,
six points that satisfy (\ref{supermario}) are as follows: the vertices $A$, $B$, and $C$, and the points diametrically opposed 
to the vertices, $A''$, $B''$, and $C''$, which all lie on the circumcircle of $\Delta ABC$.  
These six points are distinct unless $\Delta ABC$ is a right triangle, in which case, they 
reduce to four distinct points. To locate an additional point in this case that satisfies (\ref{supermario}), 
we reason as follows; suppose that $\Delta ABC$ is 
a right triangle with right angle at $C$ (thus $A=B''$ and $B=A''$).  Select
the point $D$, located on the line $BC$, with $C$ being the
midpoint of the segment $DB$. Let $C'$ be the projection
of $D$ onto $AB$.  Hence, $\Delta DCC'$ is the pedal triangle of point $D$ with 
respect to $\Delta ABC$. For simplicity of notation, if we set $BC=a$ and $AC=b$, then
$a^2+b^2=4R^2$, where $R$ is the circumradius of $\Delta ABC$.  
Clearly, $\Delta ABC$ is similar to $\Delta DBC'$.  From this, it is easy
to deduce that $DC'=\frac{ab}{R}$ and that $BC'=\frac{a^2}{R}$.  Consequently,
$OC'=\frac{a^2}{R}-R$.  Since $\Delta DOC'$ is a right triangle, $OD^2=OC'^2+DC'^2=2a^2+R^2$.
Therefore, $\frac{OD^2-R^2}{4R^2}=\frac{a^2}{2R^2}$. As for the area of 
the pedal triangle, we get
\begin{eqnarray}\label{LeBron}
|\Delta DCC'|=\frac{DC'\cdot DC\cdot\sin({\angle\!\!\!)}\,A)}{2}
=\frac{ab}{R}\cdot\frac{a}{2}\cdot\frac{a}{2R}
=\frac{a^3b}{4R^2},
\end{eqnarray}
which in combination with the fact that $|\Delta ABC| = \frac{ab}{2}$ implies
$\frac{|\Delta DCC'|}{|\Delta ABC|}=\frac{a^2}{2R^2}$. From this computation we can see 
that $D$ satisfies (\ref{supermario}).   

Now, in general, it is known that any quadratic equation 
in terms of $x$ and $y$ either has no solutions, or has as its graph a conic section. 
However, since there are points satisfying (\ref{mariokart}), the latter must be true. Thus, the shape of the locus 
of points satisfying (\ref{mariokart}) is a conic, meaning that the locus 
of points satisfying (\ref{supermario}) is either a point, two 
intersecting lines, a parabola, a hyperbola, a circle, an ellipse, or the 
whole plane (if all the lambdas are zero). Using the seven points,
$A$, $B$, $C$, $A''$, $B''$, $C''$, $O$, when $\Delta ABC$ is not right,
and the six points $A$, $B$, $C$, $D$, $C''$, $O$, when $\Delta ABC$ is right,
identified earlier as belonging to the
locus, one can eliminate all of the possible types of conics 
except for the whole plane. This means that for every $P$ in the plane
(\ref{supermario}) holds.

\section{The Area of an Antipedal Triangle}

Theorem~\ref{luigi} provides us with an efficient formula to compute 
the area of a pedal triangle given the geometry of the reference triangle 
and the location of the pedal point. This is also useful for other purposes, 
such as computing the area of an {\it antipedal} triangle in terms 
of the geometry of the reference triangle and the location of the 
antipedal point. Before proceeding with the proof of Theorem~\ref{tata5},
we prove a useful result on homotopic triangles. We do so by making use 
of a few rudiments of vector calculus.  For a reference, we refer the reader
to any multi-variable calculus textbook.

Given a triangle $A_1A_2A_3$ along with a triangle $B_1B_2B_3$ inscribed 
in it, we describe a procedure for obtaining a triangle, $C_1C_2C_3$, 
that is inscribed in $\Delta B_1B_2B_3$ and is homotopic to $\Delta A_1A_2A_3$.
By definition, two triangles are called {\it homotopic} if their sides are parallel.

\begin{proposition}\label{tata111}
Let $\Delta A_1A_2A_3$ be arbitrary and assume that $B_1\in A_2A_3$,
$B_3\in A_2A_1$, $B_2\in A_1A_3$ (see Figure~5). Take $C_1\in B_2B_3$, 
$C_2\in B_1B_3$, $C_3\in B_1B_2$ such that
\begin{eqnarray}\label{Def-Kk}
\frac{A_2B_3}{B_3A_1}=\frac{B_2C_3}{C_3B_1},\quad
\frac{A_3B_1}{B_1A_2}=\frac{B_3C_1}{C_1B_2},\quad
\frac{A_1B_2}{B_2A_3}=\frac{B_1C_2}{C_2B_3},
\end{eqnarray}
\noindent Then $\Delta A_1A_2A_3$ and $\Delta C_1C_2C_3$ are
homotopic and, in addition, $|\Delta B_1B_2B_3|$ is the geometric 
mean of $|\Delta A_1A_2A_3|$ and $|\Delta C_1C_2C_3|$, i.e.
\begin{eqnarray}\label{SS1}
|\Delta B_1B_2B_3|^2=|\Delta A_1A_2A_3|\cdot|\Delta C_1C_2C_3|. 
\end{eqnarray}
Conversely, if $\Delta A_1A_2A_3$ and $B_1\in A_2A_3$, $B_3\in A_2A_1$, 
$B_2\in A_1A_3$ are given and $C_1\in B_2B_3$, $C_2\in B_1B_3$, 
$C_3\in B_1B_2$ are such that $\Delta A_1A_2A_3$ and $\Delta C_1C_2C_3$ 
are homotopic, then (\ref{Def-Kk}) and (\ref{SS1}) hold. 
\end{proposition}

\begin{center}
\includegraphics[scale=0.70]{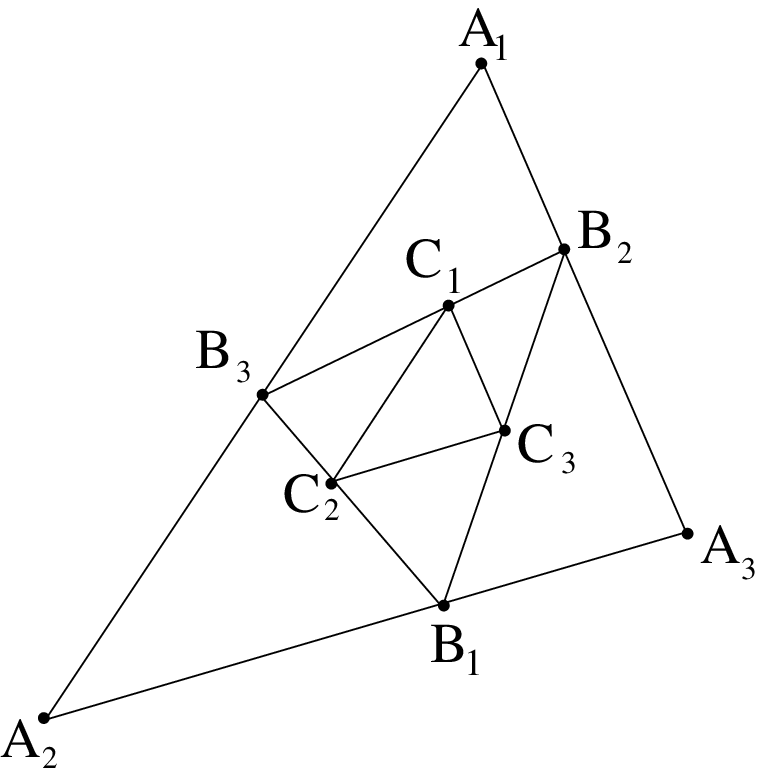}
\end{center}
\centerline{{\bf Figure~5}}
\vskip 0.15 in

\begin{proof} An affine transformation of the plane into itself
consists of a linear transformation followed by a translation and it 
has the following properties: maps lines into lines, parallel
lines into parallel lines, and preserves the ratio of line segments 
determined by points on a line.

Thus it suffices to prove Proposition~\ref{tata111} for the particular 
triangle $A_1A_2A_3$: $A_1=(0,1), A_2=(0,0), A_3=(1,0)$, since any other 
triangle can be transformed via an affine transformation into this 
particular triangle while preserving the desired properties.  
In addition, let $B_1, B_2, B_3, C_1, C_2, C_3$ be as in 
Proposition~\ref{tata111}.  We set
\begin{eqnarray}\label{Def-K}
k_1:=\frac{A_3B_1}{B_1A_2}=\frac{B_3C_1}{C_1B_2},\quad
k_2:=\frac{A_1B_2}{B_2A_3}=\frac{B_1C_2}{C_2B_3},\quad
k_3:=\frac{A_2B_3}{B_3A_1}=\frac{B_2C_3}{C_3B_1}.
\end{eqnarray}
\noindent If $M, N, P$ are three collinear points, with
coordinates $M(m_1,m_2)$, $P(p_1,p_2)$, and $N$ between $M$ and $P$, satisfying
$\frac{MN}{NP}=k$, for some real, positive constant $k$, then $N$ has 
coordinates 
\begin{eqnarray}\label{ratio}
N=\Bigl(\frac{m_1+kp_1}{1+k},\frac{m_2+kp_2}{1+k}\Bigr).
\end{eqnarray}
\noindent This fact, in combination with (\ref{Def-K}) yields 
\begin{eqnarray}\label{B}
B_1=\Bigl(\frac{1}{1+k_1},0\Bigr),\quad
B_2=\Bigl(\frac{k_2}{1+k_2},\frac{1}{1+k_2}\Bigr),\quad
B_3=\Bigl(0,\frac{k_3}{1+k_3}\Bigr).
\end{eqnarray}
\noindent Furthermore,
\begin{eqnarray}\label{C}
&& C_1=\Bigl(\frac{\frac{k_1k_2}{1+k_2}}{1+k_1},
\frac{\frac{k_3}{1+k_3}+\frac{k_1}{1+k_2}}{1+k_1}\Bigr),\quad
C_2=\Bigl(\frac{\frac{1}{1+k_1}}{1+k_2},
\frac{\frac{k_2k_3}{1+k_3}}{1+k_2}\Bigr),\quad
\nonumber\\[4pt]
&& \hskip 0.60in
C_3=\Bigl(\frac{\frac{k_2}{1+k_2}+\frac{k_3}{1+k_1}}{1+k_3},
\frac{\frac{1}{1+k_2}}{1+k_3}\Bigr).
\end{eqnarray}
\noindent It is obvious that 
\begin{eqnarray}\label{|A|}
|\Delta A_1A_2A_3|=\frac{1}{2}. 
\end{eqnarray}
Next, using vector calculus, we will compute the areas of $\Delta B_1B_2B_3$
and $\Delta C_1C_2C_3$.  More specifically, using the fact that the area of
a triangle spanned by two vectors is equal to half the norm of their cross
product, we can write,
\begin{eqnarray}\label{|B|}
|\Delta B_1B_2B_3|&=&{\textstyle{\frac{1}{2}}}
\Bigl\|\overrightarrow{B_1B_2}\times
\overrightarrow{B_1B_3}\Bigr\|\nonumber
\\[4pt]
&=&\frac{k_1k_2k_3+1}{2(1+k_1)(1+k_2)(1+k_3)}.
\end{eqnarray}
\noindent A similar reasoning applies to $\Delta C_1C_2C_3$, namely
\begin{eqnarray}\label{|C|}
|\Delta C_1C_2C_3|={\textstyle{\frac{1}{2}}}
\Bigl\|\overrightarrow{C_1C_2}\times
\overrightarrow{C_1C_3}\Bigr\|=
\frac{(k_1k_2k_3+1)^2}{2(1+k_1)^2(1+k_2)^2(1+k_3)^2}.
\end{eqnarray}
\noindent Identity (\ref{SS1}) now follows by combining (\ref{|A|}), 
(\ref{|B|}), and (\ref{|C|}), thus completing the proof of the first part 
of Proposition~\ref{tata111}.

Finally, the converse statement (as recorded in the last part of 
the proposition) follows from the uniqueness of a triangle homotopic 
with ${\triangle}A_1A_2A_3$ and inscribed in ${\triangle}B_1B_2B_3$, 
plus what we have proved so far. 
\end{proof}

\noindent{\it Proof of Theorem~\ref{tata5}.} 
If $\Delta\,DEF$ is the pedal triangle of the 
point $K_1$ with respect to $\Delta\,ABC$, then
\begin{eqnarray}\label{Saha1}
{\angle\!\!\!)} FK_1A+{\angle\!\!\!)} K_1AF=\frac{\pi}{2}.
\end{eqnarray}

\begin{center}
\includegraphics[scale=0.6]{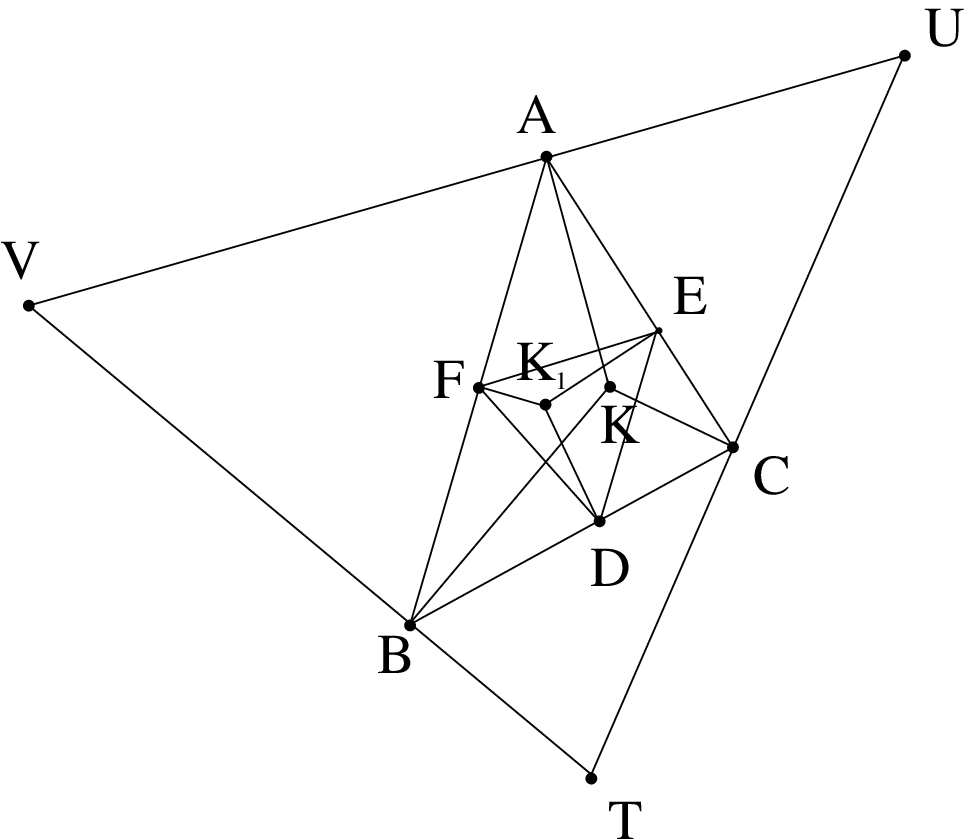}
\end{center}
\centerline{{\bf Figure~6}}
%\vskip 0.15 in

\newpage

\noindent However, because $K$ is the isogonal of $K_1$, this means that
\begin{eqnarray}\label{Saha2}
{\angle\!\!\!)} FK_1A+{\angle\!\!\!)} KAE=\frac{\pi}{2}.
\end{eqnarray}
\noindent Keeping in mind that the quadrilateral $AFK_1E$ can be inscribed
in a circle, (\ref{Saha2}) means that $AK\bot EF$, therefore $VU\|EF$.
Similar reasoning can be done to show that $VT\|DF$ and
$TU\|DE$. This implies that $\Delta TUV$ and $\Delta DEF$
are homotopic. From Proposition~\ref{tata111} we obtain
\begin{eqnarray}\label{boise4}
|\Delta DEF|\cdot|\Delta TUV|=|\Delta ABC|^2,
\end{eqnarray}
\noindent Theorem~\ref{luigi} implies
\begin{eqnarray}\label{ManU1}
\frac{|\Delta DEF|}{|\Delta ABC|}=\frac{|R^2-{OK_1}^2|}{4R^2}.
\end{eqnarray}
\noindent Therefore,  
$\frac{|\Delta TUV|}{|\Delta ABC|}=\frac{|\Delta ABC|}{|\Delta DEF|}=
\frac{4R^2}{|R^2-{OK_1}^2|}$, as claimed.

\vskip 0.1in

\vskip 0.08in
\noindent --------------------------------------
\vskip 0.10in

\noindent California Institute of Technology

\noindent MSC 700, Pasadena, CA 91126, USA

\noindent {\tt e-mail}: {\it amitrea\@@caltech.edu}

\end{document}